%

%9-12-96
\input amstex
\magnification=1200
\loadmsam
\loadmsbm
\loadeufm
\loadeusm
\UseAMSsymbols

\hsize=6.9truein
\hoffset=-0.11truein
\vsize=8.9truein
\voffset=-0.2truein

\def\leftitem#1{\item{\hbox to\parindent{\enspace#1\hfill}}}

\def\boxit#1#2{\hbox{\vrule
	\vtop{%
	\vbox{\hrule\kern#1%
	\hbox{\kern#1#2\kern#1}}%
	\kern#1\hrule}%
	\vrule}}

\def\leaderfill{\leaders\hbox to 1em{\hss.\hss}\hfill}

\parskip=\medskipamount
\document
%nopagenumbers

\input epsf
%centerline{\epsfbox{pic2-1.eps}}

%input psfig

%centerline{\bf }

\centerline{\bf Automorphisms of Thurston's Space of Measured
Laminations}

\centerline{\bf Feng Luo}

The purpose of the note is to give a characterization of the
action of the mapping class group on Thurston's space of 
measured laminations.

1. We begin with some abstract definitions. Suppose $X$ is a topological
space and $\Cal F$ is a collection of real valued (or complex valued)
functions on $X$. We say that $\Cal F$ defines an  \it
$\Cal F$-structure  \rm ($X$, $\Cal F$)
on $X$ if the topology on $X$ is the weakest topology so that each
element in $\Cal F$ is continuous, i.e., the collection
\{$f^{-1}(U) | $ $U$ open in $\bold R$, $f \in \Cal F$\} forms a subbasis for the
topology on $X$. For instance, take a smooth manifold
$X$ and let $\Cal F$ be the set of all smooth functions on $X$.
Then $(X, \Cal F)$ is the smooth structure on $X$.
An \it automorphism
\rm of a  structure ($X$, $\Cal F$) is a self-homeomorphism $\phi$ of
$X$ so that $\phi^*(\Cal F) = \Cal F$ where $\phi^*(\Cal F)
=\{ f \circ \phi | f \in \Cal F\}$.

2. Suppose now that $\Sigma$ =$\Sigma_{g, r}$ is a compact orientable
surface of genus $g$ with $r$ many boundary components ($r \geq 0$). Let 
$S(\Sigma)$ be the set of isotopy classes of homotopically non-trivial,
not boundary parallel, unoriented simple loops in $\Sigma$. Given
$\alpha$ and $\beta$ in $\Sigma$, their \it geometric intersection number, \rm
denoted by $I(\alpha, \beta)$, is the minimal number of intersection points
between their representatives, i.e.,
 $I(\alpha, \beta) =$ min$\{|a \cap b| | a \in
\alpha, b \in \beta\}$.  Thurston's space of 
(compactly supported) measured laminations on $\Sigma$, denoted by
$ML(\Sigma)$, is defined as follows. Given $\alpha$ $\in$ $S(\Sigma)$, let $
I_{\alpha}$ be
the function defined on $S(\Sigma)$ sending  $\beta$ to $I(\alpha, \beta)$. The space
$ML(\Sigma)$ is the closure of $\bold Q_{>0} \{I_{\alpha} |\alpha \in
S(\Sigma)$\} in $\bold R^{S(\Sigma)}$ under the product topology.
Thurston showed that $ML(\Sigma)$ is homeomorphic to the Euclidean space
$\bold R^{6g-6+2r}$ and the intersection pairing $I: S(\Sigma)
\times S(\Sigma) \to \bold R$ extends to a continuous pairing
$I: ML(\Sigma) \times ML(\Sigma) \to \bold R$ so that
$I(k_1 x_1, k_2 x_2) = k_1k_2 I(x_1, x_2)$ for $k_1, k_2
\in \bold R_{\geq 0}$. (See [Bo], [FLP],  [Re], [Th] and
others for a proof.) In particular, for each $\alpha$ in $S(\Sigma)$,
the map $I_{\alpha}$ from $ML(\Sigma)$ to $\bold R$ sending $m$ to 
 $I(\alpha, m)$ is continuous and the collection $\Cal F$= 
\{$I_{\alpha}$ $| \alpha \in S(\Sigma)$\} forms an $\Cal F$-structure on $ML(\Sigma)$.
According to [Th], the structure is called the \it piecewise
integral linear \rm structure on $ML(\Sigma)$. See also [Lu1].

Our result is the following.

{\bf Theorem 1.} \it  Suppose $\Sigma$ is a compact surface with or without
boundary whose Euler characteristic is negative. Then any automorphism
of the piecewise integral linear structure on the space of measured
laminations $ML(\Sigma)$ is induced by a  self-homeomorphism of the surface. \rm

3. Proof of theorem 1.

Let $\phi$ be an  automorphism of ($ML(\Sigma)$, $\Cal F$). Then $\phi$ 
induces a
bijection $\psi$ of $S(\Sigma)$ by the equation $I_{\alpha} \circ \phi = I_{\psi(\alpha)}$.

We shall first show that $\psi$ is induced  by a self-homeomorphism of the
surface. To this end, let us recall that two classes $\alpha$ and $\beta$ in
$S(\Sigma)$ are called \it disjoint, \rm denoted by $\alpha$ $\cap$ $\beta$ $= \emptyset$,
if $\alpha$ $\neq$ $\beta$ and $I(\alpha, \beta)$ = 0. By counting dimension of $I_{\alpha}^{-1}
(0)$, we shall prove that $\psi$ preserves the disjoint relation on $S(\Sigma)$. Now
by a result on the automorphism of ($S(\Sigma)$, $\cap$) (the automorphisms
of the curve complex, [Iv], [Ko], [Lu2]), we see that $\psi$ is induced
by a self-homeomorphism of the surface.

Given $\alpha$ in $S(\Sigma)$, let $Z_{\alpha} = I_{\alpha}^{-1}(0)
\subset ML(\Sigma)$. By using
the Dehn-Thurston coordinate associated to a 3-holed sphere decomposition
of the surface so that $\alpha$ is a decomposing loop ([FLP], [PH]), 
we see that the dimension  dim($Z_{\alpha}$)
is dim($ML(\Sigma)$) -1 (only the intersection
coordinate with $\alpha$ vanishes).  

{\bf Lemma 2.} \it Two elements $\alpha$ and $\beta$ in $S(\Sigma)$ are disjoint
if and only if dim($Z_{\alpha} \cap Z_{\beta}$) = dim $ML(\Sigma)$ -2. \rm

{\bf Corollary 3.} \it The bijection $\psi$ from $S(\Sigma)$ to $S(\Sigma)$ preserves the
disjointness. \rm

Indeed, the equation $I_{\alpha} \circ \phi = I_{\psi(\alpha)}$
shows that $\phi^{-1}(Z_{\alpha}) = Z_{\psi(\alpha)}$.

\noindent
\it Proof of lemma 2. \rm
We may assume that there exist disjoint elements in $S(\Sigma)$, i.e.,
dim($ML(\Sigma)$) $\geq 4$. Clearly, if $\alpha$ is disjoint from $\beta$, then
dim($Z_{\alpha} \cap Z_{\beta}$) = dim($ML(\Sigma)) -2$. This can be seen
by considering the Dehn-Thurston coordinate associated
to a 3-holed sphere decomposition so that both $\alpha$ and $\beta$
are decomposing loops. We now prove that if $\alpha \cap \beta
\neq \emptyset$, then dim($Z_{\alpha} \cap Z_{\beta})
\leq $ dim $ML(\Sigma) -3$. 

To see this, take $a \in \alpha$ and $b  \in \beta$ so that
$|a \cap b | = I(\alpha, \beta) > 0$. Let $N$ be a small regular
neighborhood of $a \cup b$. If $N$ has null homotopic boundary components
in $\Sigma$, add the disc bounded by the boundary component to $N$.
As a result, we obtain a connected subsurface $\Sigma'$ whose boundary
components are essential in $\Sigma$. Since $\alpha \cap \beta \neq \emptyset$,
the Euler characteristic of $\Sigma'$ is negative and $\Sigma'$ $\neq \Sigma_{0,3}$,
i.e., dim$(ML(\Sigma')) \geq 2$. Furthermore, $\alpha$ and
$\beta$ form a surface filling pair in $\Sigma'$, i.e., $I(\alpha, m) + I(\beta, m)
> 0$ for all $m \in ML(\Sigma')$. This implies that if $m \in ML(\Sigma)$
so that $I(m, \alpha) + I(m, \beta) = 0$, then $m$ is supported in $\Sigma -
\Sigma'$, i.e., there exist $m' \in ML(\Sigma -\Sigma')$
and some boundary components of $\alpha_1, ..., \alpha_n$ of
$\Sigma'$ so that $m$ is the disjoint union $m' \alpha_1^{k_1} ... \alpha_n
^{k_n}$ where $k_i \in \bold R_{\geq 0}$. If $\Sigma -\Sigma'$ consists
of annuli, then clearly $Z_{\alpha} \cap Z_{\beta} =\{0\}$. The result
follows. If otherwise, choose a 3-holed sphere decomposition of $\Sigma'$ and
extend it to a 3-holed sphere decomposition of $\Sigma$. For each isotopy
class  $\gamma$ of a boundary component of $\Sigma'$, we have $I(m, \gamma)
=0$ for all $m \in Z_{\alpha} \cap Z_{\beta}$. Thus, by counting the
Dehn-Thurston coordinates associated to the 3-holed sphere decomposition,
we obtain dim($Z_{\alpha} \cap Z_{\beta}) \leq$  dim($ML(\Sigma)$) -dim$(ML(
\Sigma')) -1$ $ \leq$ dim$(ML(\Sigma)) -3$.
$\square$

Now if dim$ML(\Sigma) \geq 4$ and $\Sigma \neq \Sigma_{1,2}$, then by theorem
1(a) of [Lu2] (see also [Iv], [Ko]) there exists a self-homeomorphism $f$ 
of $\Sigma$ so that $f_*^{-1}(\alpha) = \psi(\alpha)$ for all $\alpha$ in 
$S(\Sigma)$.
In particular, $I_{\alpha} \circ \phi = I_{\alpha} \circ f_*$. Since
the map from $ML(\Sigma)$ to $\bold R_{\geq 0}^{S(\Sigma)}$ sending
$m$ to $(I_{\alpha}(m))_{\alpha \in S(\Sigma)}$ is an embedding, we obtain
$\phi = f_*$.

It remains to deal with the surfaces $\Sigma$ = $\Sigma_{1,2}$, $\Sigma_{1,1}$
or $\Sigma_{0,4}$. For surface $\Sigma_{1,2}$, we shall prove that
$\psi: S(\Sigma_{1,2}) \to S(\Sigma_{1,2})$ preserves the classes
represented by separating simple loops. Assume this, then theorem 1(b)
of [Lu2] shows that $\psi$ is induced by a self-homeomorphism of the
surface. Thus, the above argument goes through.

Suppose otherwise that $\psi$ sends a separating class to a non-separating
class. We shall derive a contradiction by relating $ML(\Sigma_{1,2})$
to $ML(\Sigma_{0,5})$. Let $\tau$ be an hyper-elliptic involution
of $\Sigma_{1,2}$ with four fixed points so that the quotient space
$\Sigma_{1,2}/\tau$ is the disc $\bold D^2$ with four branch points. 
It is known by the work of Birman [Bi] and Viro [Vi] that $\tau(s)$ is isotopic
to $s$ for each simple loop $s$ not homotopic into $\partial \Sigma_{1,2}$.
In particular, we obtain $\tau_*(m) = m$ for all $m \in ML(\Sigma_{1,2})$.
Let $\pi : \Sigma_{1,2} \to \bold D^2$ be the quotient map. Consider
 $\Sigma_{0,5}$
as the disc $\bold D^2$ with a regular neighborhood  $N(B)$ of the the
branched point set $B$ removed, i.e., $\Sigma_{0,5} = \bold D - int(N(B))$.
Define $p : S(\Sigma_{0,5}) \to ML(\Sigma_{1,2})$ by sending the
isotopy class $[a]$ to the measured lamination $[\pi^{-1}(a)]$. Note
that if $\pi^{-1}(a)$ is connected, then it is a  separating loop
and if $\pi^{-1}(a)$ is not connected, then it
is a union of two parallel copies of a non-separating simple
loops.  This map
$p$ extends to a homeomorphism, still denoted by $p$, from $ML(\Sigma_{0,5})$
to $ML(\Sigma_{1,2})$ so that $I(p(m_1), p(m_2)) = 2 I(m_1, m_2)$ for
all $m_1, m_2 \in ML(\Sigma_{0,5})$.

Now consider the homeomorphism $\phi': ML(\Sigma_{0,5}) \to ML(\Sigma_{0,5})$
given by $p^{-1} \phi p$. Since $I_{\alpha} \phi = I_{\psi(\alpha)}$
for all $\alpha$ $\in S(\Sigma_{1,2})$, we obtain $\lambda I_{\alpha} \circ
\phi'
= I_{\psi'(\alpha)}$ for all $\alpha \in S(\Sigma_{0,5})$ where
$\psi': S(\Sigma_{0,5}) \to S(\Sigma_{0,5})$ is a bijection and
$\lambda = 1$ or  $1/2$ or $2$ depending on the components of 
$\pi^{-1}(\alpha)$ and $\pi^{-1}(\phi'(\alpha))$
 being separating or not. By the assumption that $\psi$ sends some
non-separating simple loops to separating  ones, the function 
$\lambda$ is not a constant. Due to the  equation 
$\lambda I_{\alpha} \circ \phi'
= I_{\psi'(\alpha)}$, the map $\phi'$  preserves the set \{$Z_{\alpha}$
$| \alpha \in S(\Sigma_{0,5})$\}.
By lemma 2, we see that $\psi'$ preserves the disjointness. Thus
$\psi'$ is induced by a self-homeomorphism $h$ of $\Sigma_{0,5}$. In
particular we obtain $\lambda I_{\alpha} \circ \phi' = I_{\alpha}  \circ h$ for
all $\alpha$. Since the set of rational multiples of $S(\Sigma_{0,5})$
is dense in $ML(\Sigma_{0,5})$ and both $\phi'$ and $h$ are
are homogeneous, it follows that $\phi' = k h$ for some
fixed constant $k=1$ or $1/2$ or $2$. This contradicts the assumption that
$\lambda$ is not a constant.

Finally, we show that any automorphism of $(ML(\Sigma_{1,1}), \Cal F)$
and $(ML(\Sigma_{0,4}), \Cal F)$ is induced by a surface homeomorphism.
Since the structures $(ML(\Sigma_{1,1}), \Cal F)$ and $(ML(\Sigma_{0,4}),
\Cal F)$ are isomorphic, we shall deal with the case $\Sigma_{1,1}$ only.

Let us first identify both $ML(\Sigma_{1,1})$ and $S(\Sigma_{1,1})$
with the first homology groups. Let $i : S(\Sigma_{1,1})
\to H_1(\Sigma_{1,1}, \bold Z)/\pm 1$ be the natural map sending
an isotopy class to the corresponding homology classes. It is well
known that the map is a bijection from $S(\Sigma_{1,1})$ to
$\Cal P/ \pm 1$ where $\Cal P$ is the set of primitive elements
in $H_1(\Sigma_{1,1}, \bold Z) \cong \bold Z^2$. Furthermore, by taking a
$\bold Z$-basis for $H_1(\Sigma_{1,1}, \bold Z)$, each $i(\alpha)$
can be written as $\pm (a, b)$ where $a, b$
are relatively prime integers. Under this identification, the
intersection number $I(\alpha_1, \alpha_2) = |a_1 b_2 - a_2 b_1|$
where $\alpha_i = \pm (a_i ,b_i)$.
In particular, this shows that $ML(\Sigma_{1,1})$ can be naturally
identified with $H_1(\Sigma_{1,1}, \bold R)/\pm 1 \cong \bold R^2/\pm 1$
so that the
above intersection number formula still holds.

The action of self-homeomorphisms on $ML(\Sigma_{1,1})$ is
induced by the $GL(2, \bold Z)$ action on $\bold R^2/\pm 1$.  Thus, it
remains to show that if $\phi =(\phi_1, \phi_2)
:\bold R^2/\pm 1 \to \bold R^2/\pm 1$ is a self-homeomorphism so that
for each pair of relative prime integers $(a, b) \in \Cal P$ there exists
a new pair $(a', b') \in \Cal P$ satisfying
$| a \phi_1(x,y) - b \phi(x, y) | = |a'x - b'y|$ for all $(x,y) \in \bold R^2$,
then $\phi$ is induced by an element in $GL(2, \bold Z)$. By taking
$(a,b)$ to be $(1,0)$ and $(0,1)$, we see that $|\phi_1(x,y)| = | a_1x + b_1y|$
and $|\phi_2(x,y)| = |a_2x + b_2y|$. Since $\phi$ is a homeomorphism,
$a_1b_2 - a_2b_1 \neq 0$. The goal is to show that $a_1b_2 -a_2 b_1 =
\pm 1$. Since both $\phi_1$, $\phi_2$ are continuous and
$|\phi_1(x,y) \pm \phi_2(x,y)|$ is of the form $|ax + by|$, 
it follows that $\phi_i(x,y) = \pm (a_ix + b_iy)$ for $i=1,2$.
Now if $|a_1 b_2 - a_2 b_1| \geq 2$, then one can find $(a, b) \in \Cal P$
so that $a\phi_1 - b\phi_2$ is of the form $|cx + dy|$ where $c$ and
$d$ have a common non-trivial divisor. This contradicts the assumption.
$\square$

4. One consequence of the proof of the theorem 1 is the following
characterization of the action of the mapping class group on the
projectivized measured lamination space $PML(\Sigma) = ML(\Sigma) -\{0\}/
\bold R_{>0}$.

{\bf Theorem 4.} (Automorphisms of the projective measured lamination spaces)
\it Suppose $\Sigma$ is a compact orientable surface so that
dim$(ML(\Sigma)) \geq 2$ and $\Sigma \neq \Sigma_{1,2}$. For each
$\alpha \in $ $S(\Sigma)$, let $P_{\alpha}$ be the image of
$\{ m \in ML(\Sigma) -\{0\}| I(m, \alpha) = 0\}$ in $PML(\Sigma)$. If
$\phi$ is a self-homeomorphism of the projective measured lamination
space $PML(\Sigma)$ preserving the collection $\{ P_{\alpha}
| \alpha \in S(\Sigma)\}$, then $\phi$ is induced by a self-homeomorphism
of the surface. \rm

\noindent
\it Proof. \rm By lemma 2 and the result on the automorphism of the
curve complex, we see that there exists a self-homeomorphism $f$ of
the surface so that $f_* \phi^{-1}: PML(\Sigma) \to PML(\Sigma)$ sends
each $P_{\alpha}$ to $P_{\alpha}$. The image of $\Cal P(\alpha)$
of $\alpha \in S(\Sigma)$ in $PML(\Sigma)$ can be expressed 
as a finite intersection $P_{\alpha_1} \cap P_{\alpha_2} \cap ... \cap 
P_{\alpha_k}$. Thus
$f_*\phi^{-1}$ is the identity map on the set \{$\Cal P(\alpha)$
$| \alpha \in S(\Sigma)\}$. Since the set \{ $\Cal P(\alpha)$ 
$| \alpha \in S(\Sigma)\}$ is dense in 
$PML(\Sigma)$, it follows that $\phi = f_*$.
$\square$

5. \it Remark. \rm  The theorem is  valid  for $\Sigma_{1,2}$
if we assume that the self-homeomorphism $\phi$ preserves the
subset $\{P_{\alpha} | \alpha $ is a separating class \}. Otherwise,
it is false. See [Lu2].

6. Similar automorphism results hold for the Teichmuller
space and $SL(2, \bold R)$ characters. For simplicity, we state the result for
the Teichmuller space. The proof is essentially the same as above
and will be omitted.  
Let $T(\Sigma)$  be the space of isotopy classes of hyperbolic
metrics with cusp ends on $int(\Sigma)$. For each $\alpha$ $\in S(\Sigma)$,
let $l_{\alpha}: T(\Sigma) \to \bold R$ be the geodesic length
function sending a metric $m$ to the length of $m$-geodesic in
$\alpha$. The work of Fricke-Klein [FK] shows that the collection
\{$ l_{\alpha} | \alpha \in S(\Sigma)\}$ forms an $\Cal F$-structure
on the Teichmuller space.

{\bf Theorem 5.} \it Suppose $\Sigma$ is a compact surface of negative
Euler characteristic. Then any automorphism of $(T(\Sigma), \Cal F)$
is induced by a self-homeomorphism of the surface. \rm

The key step in the proof is to show a result similar to lemma 2.
In this case, it is the Margulis lemma that $\alpha \cap \beta
=\emptyset$ if and only if $inf \{ l_{\alpha} + l_{\beta} \} = 0$
on $T(\Sigma)$.

7. Currently, we are unable to solve the automorphism problem
for the variety of characters of $SL(2,\bold C)$ representations of
 a closed surface group with respect to the structure of the trace
functions $\{ tr_{\alpha} | \alpha \in S(\Sigma)\}$. 
 Here
$tr_{\alpha}$ sends a character $\chi$ to $\chi(\alpha)$.
See [CS]
for an introduction to the subject. 
The main difficulty is due to the lacking  of intrinsic characterization
of disjointness $\alpha \cap \beta = \emptyset$ in terms of the
trace functions $tr_{\alpha}$ and $tr_{\beta}$.

8. \it Acknowledgement.  \rm The work is partially supported by the NSF.

\centerline{\bf Reference}

[Bi] Birman, J.S.,  Braids, links, and mapping class groups.
Ann. of Math. Stud., 82, Princeton Univ. Press, Princeton, NJ, 1975

[Bo] Bonahon, F.: The geometry  of Teichm\"uller space via geodesic
currents. Invent. Math. {\bf 92} (1988), 139-162

[CS] Culler, M., Shalen, P.: Varieties of group representations and
splittings of 3-manifolds. Ann. Math. {\bf 117} (1983), 109-146

[FK] Fricke, R., Klein, F.: Vorlesungen  \"uber die Theorie der
Automorphen Functionen. Teubner, Leipizig, 1897-1912

[FLP] Fathi, A., Laudenbach, F., Poenaru, V.: Travaux de Thurston sur les
surfaces. Ast\'erisque {\bf 66-67}, Soci\'et\'e Math\'ematique de France, 1979

[Iv] Ivanov, N.V.: Automorphisms of complexes of curves and of Teichm\"uller
 spaces. Inter. Math. Res. Notice, 14, (1997), 651-666

[Ko] Korkmaz, M: Complexes of curves and the mapping class groups, Top.
 \& Appl., to appear.

[Lu1] Luo, F.: Simple loops on surfaces and their intersection numbers,
preprint, 1997. The Research announcement appeared in  Math. Res. Letters,
{\bf 5}, (1998), 47-56.

[Lu2] Luo, F.:  Automorphisms of the curve complex, Topology,
to appear.

[PH] Penner, R.,  Harer, J.: Combinatorics of train tracks.
Annals of Mathematics Studies, 125. Princeton University Press, Princeton,
NJ, 1992.

[Re] Rees, M.: An alternative approach to the ergodic theory of
measured foliations on surfaces. Ergodic Theory Dynamical Systems 1 (1981),
no. 4, 461-488.

[Th] Thurston, W.: On the geometry and dynamics of diffeomorphisms of
surfaces. Bul. Amer. Math. Soc. {\bf 19} (2) (1988), 417-438       

[Vi] Viro, O.: Links, two-fold branched coverings and braids,
 Soviet Math. Sbornik, 87 no 2 (1972), 216-228

Department of Mathematics

Rutgers University

Piscataway, NJ 08854

\end

\end